\documentclass[12pt,psamsfonts,leqno]{amsart}
\usepackage{amssymb,eucal,xy}
\xyoption{all}
\xyoption{ps}

\theoremstyle{plain}

\newtheorem{lemma}{Lemma}[section]
\newtheorem{theorem}{Theorem}
\newtheorem{proposition}[lemma]{Proposition}
\newtheorem{corollary}[lemma]{Corollary}
\newtheorem{claim}{Claim}

\newtheorem*{stat}{\name}
\newcommand{\name}{testing}

\theoremstyle{definition}
\newtheorem{definition}[lemma]{Definition}

\newtheorem{problem}{Problem}

\theoremstyle{remark}
\newtheorem{remark}[lemma]{Remark}

\newenvironment{all}[1]{\renewcommand{\name}{#1}\begin{stat}}
                        {\end{stat}}

\newcommand{\qedc}{{\qed}~{\rm Claim~{\theclaim}.}}

\newenvironment{cproof}
{\begin{proof}[Proof of Claim.]}
{\qedc\renewcommand{\qed}{}\end{proof}}

\numberwithin{equation}{section}

\newcommand{\cjh}{com\-plete join-ho\-mo\-mor\-phism}
\newcommand{\cmh}{com\-plete meet-ho\-mo\-mor\-phism}
\newcommand{\kl}{{\kappa,\lambda}}
\newcommand{\ab}{{\vec\xa,\vec\xb}}
\newcommand{\klip}{$\seq{\kl}$-in\-ter\-po\-la\-tion
property}
\newcommand{\pscp}{pse\-udo-com\-ple\-men\-ted}

\newcommand{\cb}{co-Brouwerian}
\newcommand{\ccb}{conditionally co-Brou\-wer\-ian}
\newcommand{\ckcb}{conditionally $\kappa$-co-Brou\-wer\-ian}

\newcommand{\Ccb}{Conditionally co-Brou\-wer\-ian}

\newcommand{\Cj}{\mathbf{C}_{\vee}}
\newcommand{\Cm}{\mathbf{C}_{\wedge}}
\newcommand{\KK}{\mathcal{K}}

\newcommand{\eps}{\varepsilon}
\newcommand{\es}{\varnothing}
\newcommand{\into}{\hookrightarrow}
\newcommand{\onto}{\twoheadrightarrow}
\newcommand{\cep}{congruence extension property}
\newcommand{\seq}[1]{\langle#1\rangle}

\newcommand{\set}[1]{\{#1\}}
\newcommand{\setm}[2]{\set{#1\mid#2}}

\newcommand{\ol}[1]{\overline{#1}}
\newcommand{\oll}[1]{\,\overline{\!#1}}

\newcommand{\zero}{\mathbf{0}}
\newcommand{\one}{\mathbf{1}}
\newcommand{\two}{\mathbf{2}}
\newcommand{\go}{\omega}

\newcommand{\sd}{\smallsetminus}
\newcommand{\dnw}{\mathbin{\downarrow}}

\DeclareMathOperator{\Id}{Id}

\DeclareMathOperator{\Int}{Int}
\DeclareMathOperator{\Ob}{Ob}
\DeclareMathOperator{\rng}{rng}

\newcommand{\II}{\mathcal{I}}
\newcommand{\FF}{\mathcal{F}}

\newcommand{\PL}{\mathbf{PL}}
\newcommand{\LL}{\mathbf{L}}
\newcommand{\CC}{\mathbf{C}}
\newcommand{\DD}{\mathcal{D}}
\newcommand{\SL}{\mathbf{S}}
\newcommand{\MPL}[1]{$#1$-meas\-ured par\-tial lat\-tice}
\newcommand{\ML}[1]{$#1$-meas\-ured lat\-tice}

\def\cls(#1,#2){{{#1}{}_{/{#2}}}}

\DeclareMathOperator{\Con}{Con}
\DeclareMathOperator{\Res}{Res}
\DeclareMathOperator{\Conc}{Con_c}

\newcommand{\id}{\mathrm{id}}
\newcommand{\jz}{$\langle\vee,0\rangle$}
\newcommand{\jzh}{\jz-ho\-mo\-mor\-phism}
\newcommand{\jzs}{\jz-sem\-i\-lat\-tice}
\newcommand{\res}{\mathbin{\restriction}}

\newcommand{\fine}[1]{[#1]_*^{<\go}}

\DeclareMathOperator{\Fg}{F}

\newcommand{\FL}{\Fg_{\mathbf{L}}}

\newcommand{\xa}{{\boldsymbol{a}}}
\newcommand{\xb}{{\boldsymbol{b}}}
\newcommand{\xc}{{\boldsymbol{c}}}
\newcommand{\xd}{{\boldsymbol{d}}}

\newcommand{\xx}{{\boldsymbol{x}}}

\newcommand{\xC}{{\boldsymbol{C}}}
\newcommand{\xX}{{\boldsymbol{X}}}

\begin{document}

\title[Two-dimensional congruence amalgamation of lattices]%
{Join-semilattices with two-dimensional congruence amalgamation}

\author[F.~Wehrung]{Friedrich Wehrung}
\address{CNRS, FRE 2271\\
D\'epartement de Math\'ematiques\\
Universit\'e de Caen\\
14032 Caen Cedex\\
France}
\email{wehrung@math.unicaen.fr}
\urladdr{http://www.math.unicaen.fr/\~{}wehrung}

\date{\today}


\begin{abstract}
We say that a \jzs\ $S$ is \emph{\ccb}, if (1) for all nonempty subsets
$X$ and $Y$ of $S$ such that $X\leq Y$ (\emph{i.e.}, $x\leq y$ for all
$\seq{x,y}\in X\times Y$), there exists $z\in S$ such that
$X\leq z\leq Y$, and (2) for every subset $Z$ of $S$ and all $a$,
$b\in S$, if $a\leq b\vee z$ for all $z\in Z$, then there exists $c\in S$
such that $a\leq b\vee c$ and $c\leq Z$. By restricting this definition to
subsets $X$, $Y$, and $Z$ of less than $\kappa$ elements, for an infinite
cardinal $\kappa$, we obtain the definition of a \emph{\ckcb} \jzs.

We prove that for every \ccb\ lattice $S$ and every partial lattice $P$,
every \jzh\ $\varphi\colon\Conc P\to S$ can be lifted to a lattice
homomorphism $f\colon P\to L$, for some relatively complemented lattice
$L$. Here, $\Conc P$ denotes the \jzs\ of compact congruences of $P$.

We also prove a two-dimensional version of this result, and we establish
partial converses of our results and various of their consequences in
terms of congruence lattice representation problems. Among these
consequences, for every infinite regular cardinal $\kappa$ and every
\ckcb\ $S$ of size $\kappa$, there exists a relatively complemented
lattice $L$ with zero such that $\Conc L\cong S$.
\end{abstract}

\maketitle

\section{Introduction}

The present paper deals essentially with two categories of structures. The
first one is the category $\PL$ of all \emph{partial lattices} (see
Definition \ref{D:PartLatt}) and their homomorphisms (see Definition
\ref{D:HomPL}), while the second one is the category $\SL$ of all \jzs s
and \jzh s. These categories are related by the functor
$\Conc\colon\PL\to\SL$. For a partial lattice $P$, $\Conc P$ is the \jzs\
of compact congruences of $P$, see Section~\ref{S:PartLatt}.

In the last few years some effort has been put on the investigation of
the effect of the $\Conc$ functor not only on the \emph{objects} of $\PL$,
but also on the \emph{diagrams} of $\PL$, in fact, essentially on
the diagrams of the full subcategory $\LL$ of $\PL$ whose objects are all
\emph{lattices}. A complete account of the pre-1998 stages of this
research is presented in \cite{GrScC}.
Formally, a diagram of $\PL$ is a functor from a category
$\CC$ to $\PL$. Most of the results of the last years in this topic can
then be conveniently formulated \emph{via} the following definition.

\begin{definition}\label{D:FactLift}
Let $\DD$ be a diagram of partial lattices. We denote the
composition $\Conc\circ\DD$ by $\Conc\DD$. For a \jzs\ $S$ and a
partial lattice $P$, we say that a homomorphism
$\varphi\colon\Conc\DD\to S$ can be
\begin{enumerate}
\item \emph{factored through} $P$, if there are a homomorphism
$f\colon\DD\to P$ and a \jzh\ $\psi\colon\Conc P\to S$ such that
$\varphi=\psi\circ\Conc f$;

\item \emph{lifted through} $P$, if there are a homomorphism
$f\colon\DD\to P$ and an isomorphism $\psi\colon\Conc P\to S$ such that
$\varphi=\psi\circ\Conc f$.
\end{enumerate}

In (i) (resp., (ii)) above, we say that $\varphi$ can be \emph{factored
to} (resp., \emph{lifted to})~$f$.
\end{definition}

Homomorphisms between diagrams have to be understood in the categorical
sense, \emph{e.g.}, if $\DD\colon\CC\to\PL$ is a diagram of partial
lattices and if $P$ is a partial lattice, a homomorphism $f\colon\DD\to P$
consists of a family $(f_X)_{X\in\Ob\CC}$ of homomorphisms
$f_X\colon\DD(X)\to P$, for any object $X$ of $\CC$, such that if
$u\colon X\to Y$ is a morphism in~$\CC$, then $f_X=f_Y\circ\DD(u)$. Of
particular interest to us will be the case where $\DD$ consists exactly
of one partial lattice, \emph{i.e.}, $\CC$ is the trivial category with
one object and one morphism, and the case where $\DD$ is a
\emph{truncated square}, \emph{i.e.}, $\CC$ consists of distinct
objects $0$, $1$, and $2$ together with nontrivial morphisms
$e_1\colon 0\to 1$ and $e_2\colon 0\to 2$. In that case $\DD$ can be
described by partial lattices $P_0$, $P_1$, and~$P_2$, together with
homomorphisms
$f_1\colon P_0\to P_1$ and $f_2\colon P_0\to P_2$. Moreover, if $P$ is a
partial lattice, a homomorphism from $\DD$ to $P$ can then be
described by homomorphisms of partial lattices $g_i\colon P_i\to P$, for
$i<3$, such that $g_1\circ f_1=g_2\circ f_2=g_0$. The situation can be
described by the following commutative diagrams:
 \[
{
\def\labelstyle{\displaystyle}
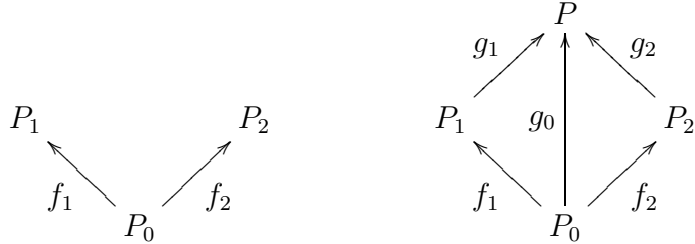
\xymatrix{
& & & & & P &\\
P_1 & & P_2 & & P_1\ar[ur]^{g_1} & & P_2\ar[ul]_{g_2}\\
& P_0\ar[ul]^{f_1}\ar[ur]_{f_2} & &
\save+<0ex,-5ex>\drop{\text{Illustrating $\DD$
and a homomorphism from $\DD$ to $P$}}\restore & &
P_0\ar[ul]^{f_1}\ar[ur]_{f_2}\ar[uu]^{g_0} &
}
}
 \]
We shall call $P_0$ (resp., $P_1$ and
$P_2$) the \emph{bottom} (resp., the \emph{sides}) of $\DD$.

Then a typical lifting result of the $\Conc$ functor is the following,
see J.~T\r{u}ma \cite{Tuma} and G.~Gr\"atzer, H.~Lakser, and
F.~Wehrung~\cite{GLW}:

\begin{theorem}\label{T:Tuma}
Let $\DD$ be a truncated square of lattices, let $S$ be a finite
distributive \jzs. Then every homomorphism $\varphi\colon\Conc\DD\to S$
can be lifted through a relatively complemented lattice.
\end{theorem}

For infinite $S$, completely different methods yield the following
result, see Theorem~C in F.~Wehrung \cite{Wehr}.

\begin{theorem}\label{T:1DimLift}
Let $K$ be a lattice, let $S$ be a distributive \emph{lattice} with
zero. Then every \jzh\ $\varphi\colon\Conc K\to S$ can be lifted.
Furthermore, a lift
$f\colon K\to L$ can be found in such a way that the following assertions
hold:
\begin{enumerate}
\item $L$ is relatively complemented.

\item The range of $f$ generates $L$ as an ideal (resp., a filter).

\item If the range of $\varphi$ is cofinal in $S$, then the range of
$f$ generates $L$ as a convex sublattice.
\end{enumerate}
\end{theorem}

We will express the repetition of the conditions (i)--(iii) above by
saying that $\varphi$ \emph{has a good lift}, although, strictly
speaking, one would need to define good ideal lifts and good filter lifts.
Moreover, this condition turns out to be somehow looser than it appears,
as, for example, it can be strengthened by many additional properties of
$L$, such as the ones listed in the statement of Proposition~20.8 of
\cite{Wehr}. For example, $L$ has \emph{definable principal congruences}.

One can then say that Theorem~\ref{T:Tuma} is a two-dimensional
lifting result for finite distributive \jzs s, while
Theorem~\ref{T:1DimLift} is a one-dimensional lifting result for
arbitrary distributive lattices with zero. In fact,
the following stronger, ``two-dimensional'' result holds, see
\cite[Theorem~D]{Wehr}:

\begin{theorem}\label{T:<2DimLift}
Let $\DD$ be a truncated square of lattices with \emph{finite} bottom, let
$S$ be a distributive lattice with zero, let
$\varphi\colon\Conc\DD\to S$ be a homomorphism. Then $\varphi$ has a
good lift.
\end{theorem}

In the `good lift' statement, the range of $f\colon\DD\to L$ has to be
understood as the union of the ranges of the images under $f$ of the
individual objects in $\DD$, while the range of $\varphi$ is the
\jzs\ generated by the union of the ranges of the images under
$\varphi$ of the individual objects in $\Conc\DD$.

As we shall see in the present paper, the statement of Theorem~3 does
\emph{not} extend to the case where the bottom of $\DD$ is an
\emph{infinite} lattice. However, we shall introduce a class of
distributive lattices with zero, the so-called \emph{\ccb} ones (see
Definition~\ref{D:CB}), that includes all finite distributive lattices.
Moreover, every complete sublattice of a complete Boolean lattice is
\ccb. For those lattices, the stronger statement remains valid, and much
more:

\begin{theorem}\label{T:1Main}
Let $P$ be a partial lattice, let~$S$ be a \ccb\ lattice. Then every
homomorphism $\varphi\colon\Conc P\to S$ has a good lift.
\end{theorem}

Now the two-dimensional version of Theorem~\ref{T:1Main}:

\begin{theorem}\label{T:2Main}
Let $\DD$ be a truncated square of partial lattices with bottom a
lattice, let~$S$ be a \ccb\ lattice. Then every homomorphism
$\varphi\colon\Conc\DD\to S$ has a good lift.
\end{theorem}

As we shall prove in Sections \ref{S:Exmples} and \ref{S:Cpl}, some of
the assumptions on $S$ are also necessary for the statements of
Theorems \ref{T:1Main} and \ref{T:2Main} to hold.

All these results imply the following corollaries:

\begin{all}{Corollary \ref{C:MoreRepr}}
Let $S$ be a distributive \jzs\ that can be expressed as a \jz-direct
limit of at most $\aleph_1$ \ccb\ lattices. Then there exists a
relatively complemented lattice $L$ with zero such that $\Conc L\cong S$.
Furthermore, if $S$ is bounded, then $L$ can be taken bounded as well.
\end{all}

This result extends a well-known result of A. Huhn \cite{Huhn89a,Huhn89b}
that states that every distributive \jzs\ of size at most $\aleph_1$ is
isomorphic to $\Conc L$ for some lattice $L$.

Our next corollary also implies a positive solution for Problem~4
of~\cite{GLW}.

\begin{all}{Corollary \ref{C:CPE}}
Let $K$ be a lattice that can be expressed as a direct union of
\emph{countably} many lattices whose congruence semilattices are \ccb.
Then $K$ embeds congruence-preservingly into some relatively
complemented lattice~$L$, which it generates as a convex sublattice.
\end{all}

We also establish \emph{relativizations} of the methods leading to
Theorems \ref{T:1Main} and \ref{T:2Main}. These statements involve a
relativized version, for every infinite cardinal $\kappa$, of the notion
of a \ccb\ lattice. We call \emph{\ckcb\ \jzs s} the resulting
objects, see Definition~\ref{D:kCB}.

\begin{theorem}\label{T:kappaRepr}
Let $\kappa$ be an infinite cardinal, let $S$ be a \ckcb\
\jzs\ of size $\kappa$. Then there exists a relatively complemented
lattice~$L$ with zero such that $\Conc L\cong S$. Furthermore, if $S$
is bounded, then $L$ can be taken bounded as well.
\end{theorem}

\section{Algebraic lattices}\label{S:AlgLatt}

We recall that in a lattice $L$, an element $a$ is
\emph{compact}, if for any nonempty upward directed subset $X$ of $L$,
$a\leq\bigvee X$ implies that $a\leq x$ for some $x\in X$. We denote by
$\KK(L)$ the join-semilattice of compact elements of~$L$. We say that
$L$ is \emph{algebraic}, if $L$ is complete and every element of $L$ is a
join of compact elements. If $L$ is an algebraic lattice, then $\KK(L)$
is a \jzs, while for every \jzs\ $S$, the lattice $\Id S$ of all ideals
of $S$ is an algebraic lattice. These transformations can be extended to
functors in a canonical way. The relevant definitions for the morphisms
are the following. For \jzs s, they are the \jzh s, while for algebraic
lattices, they are the compactness preserving \cjh s; by definition, for
complete lattices $A$ and $B$, a map $f\colon A\to B$ is a \emph{\cjh},
if $f\left(\bigvee X\right)=\bigvee f[X]$ for any subset $X$ of $A$,
while $f$ is \emph{compactness preserving}, if
$f[\KK(A)]\subseteq\KK(B)$. Then the aforementioned category equivalence can
be stated in the following condensed form:

\begin{proposition}\label{P:AlgLattSem}
The functors $S\mapsto\Id S$ and $A\mapsto\KK(A)$ define a category
equivalence between \jzs s with \jzh s and algebraic lattices with
compactness preserving \cjh s.
\end{proposition}

\section{Partial lattices}\label{S:PartLatt}

Our notations and definitions are the same as in \cite{Wehr}. If $X$ is a
subset of a \emph{quasi-ordered} set $P$ and if $a\in P$, let $a=\sup X$
(resp., $a=\inf X$) be the statement that $a$ is a majorant (resp.,
minorant) of $X$ and that every majorant (resp., minorant) $x$ of $X$
satisfies that $a\leq x$ (resp., $x\leq a$). We observe that this
statement determines $a$ only up to equivalence.

\begin{definition}\label{D:PartLatt}
A \emph{partial prelattice} is a structure
$\seq{P,\leq,\bigvee,\bigwedge}$, where $P$ is a nonempty set, $\leq$
is a quasi-ordering on~$P$, and $\bigvee$, $\bigwedge$ are partial
functions from the set $\fine P$ of all nonempty finite subsets of $P$
to $P$ satisfying the following properties:
\begin{enumerate}
\item $a=\bigvee X$ implies that $a=\sup X$,
for all $a\in P$ and all $X\in\fine P$.

\item $a=\bigwedge X$ implies that $a=\inf X$,
for all $a\in P$ and all $X\in\fine P$.
\end{enumerate}

We say that $P$ is a \emph{partial lattice}, if $\leq$ is antisymmetric.

A \emph{congruence} of~$P$ is a quasi-ordering $\preceq$ of~$P$ containing
$\leq$ such that\linebreak $\seq{P,\preceq,\bigvee,\bigwedge}$ is a partial
prelattice.
\end{definition}

For a partial lattice $P$, a congruence $\xc$ of $P$, and elements $x$,
$y$ of $P$, we shall often write $x\leq_\xc y$ instead of
$\seq{x,y}\in\xc$, and $x\equiv_\xc y$ instead of the conjunction of
$x\leq_\xc y$ and $y\leq_\xc x$. The \emph{quotient} $P/{\xc}$ has
underlying set $P/{\equiv_\xc}$, we endow it with the quotient
quasi-ordering ${\leq_\xc}/{\equiv_\xc}$ and the partial join defined by
the rule
 \begin{multline*}
 \xa=\bigvee\xX\qquad\text{if{f}}\qquad\text{there are }X\in\fine{P}
 \text{ and }a\in P\\
 \text{with }a=\bigvee X,\ \xa=a_{/{\xc}},\text{ and }
 \xX=X_{/{\xc}},
 \end{multline*}
where $a_{/{\xc}}$ denotes the equivalence class of $a$ modulo $\xc$ and
we put $X_{/{\xc}}=\setm{x_{/{\xc}}}{x\in X}$. The partial meet on
$P/{\xc}$ is defined dually.

For $a$, $b\in P$, we denote by $\Theta_P^+(a,b)$ the least congruence
$\xc$ of~$P$ such that $a\leq_\xc b$, and we put
$\Theta_P(a,b)=\Theta_P^+(a,b)\vee\Theta_P^+(b,a)$, the least congruence
$\xc$ of~$P$ such that $a\equiv_\xc b$. Of course, the congruences of
the form $\Theta_P^+(a,b)$ are generators of the join-semilattice $\Conc P$.

We shall naturally identify \emph{lattices} with partial lattices $P$ such
that $\bigvee$ and $\bigwedge$ are defined everywhere on $\fine P$.

\begin{proposition}\label{P:ConP}
Let $P$ be a partial prelattice. Then the set $\Con P$ of all congruences of
$P$ is a closure system in the powerset lattice of~$P\times P$, closed under
directed unions. In particular, it is an algebraic lattice.
\end{proposition}

We denote by $\Conc P$ the \jzs\ of all \emph{compact}
congruences of~$P$, by $\zero_P$ the least congruence of~$P$
(that is, $\zero_P$ is the quasi-ordering of~$P$), and by $\one_P$ the
largest (coarse) congruence of~$P$.

If $P$ is a lattice, then $\Con P$ is distributive, but
this may not hold for a general partial lattice $P$.

Many \jzh s will be constructed by using the following notion of
\emph{measure}.

\begin{definition}\label{D:meas}
Let $P$ be a partial lattice, let $S$ be a \jzs. A \emph{$S$-valued
measure} on $P$ is a map $\mu\colon P\times P\to S$ that satisfies the
following properties (we will write from now on $\mu(x,y)$ instead of
$\mu(\seq{x,y})$):
\begin{enumerate}
\item $\mu(x,y)=0$, for all $x$, $y\in P$ such that $x\leq y$.

\item $\mu(x,z)\leq\mu(x,y)\vee\mu(y,z)$, for all $x$,
$y$, $z\in P$.

\item $\mu(a,b)=\bigvee_{x\in X}\mu(x,b)$, for all $a$, $b\in P$
and all $X\in\fine P$ such that $a=\bigvee X$.

\item $\mu(a,b)=\bigvee_{y\in Y}\mu(a,y)$, for all $a$, $b\in P$
and all $Y\in\fine P$ such that $b=\bigwedge Y$.
\end{enumerate}
\end{definition}

We omit the easy proof of the following lemma, see also Proposition~13.1
in~\cite{Wehr}. This lemma states that the notion of measure on $P$ and the
notion of \jzh\ from $\Conc P$ are essentially equivalent.

\begin{lemma}\label{L:meas}
Let $P$ be a partial lattice, let $S$ be a \jzs. Then the following
assertions hold:
\begin{enumerate}
\item For every \jzh\ $\ol{\mu}\colon\Conc P\to S$, the map\linebreak
$\mu\colon P\times P\to S$, $\seq{x,y}\mapsto\ol{\mu}\Theta_P^+(x,y)$ is
a $S$-valued measure on $P$.

\item For any $S$-valued measure $\mu$ on $P$, there exists a unique
\jzh\
$\ol{\mu}\colon\Conc P\to S$ such that
$\mu(x,y)=\ol{\mu}\Theta_P^+(x,y)$, for all $x$, $y\in P$.
\end{enumerate}
\end{lemma}

The homomorphism $\ol{\mu}$ (the ``integral'' with respect to $\mu$)
is of course defined by the formula
 \[
 \ol{\mu}\left(\bigvee_{i<n}\Theta_P^+(x_i,y_i)\right)
 =\bigvee_{i<n}\mu(x_i,y_i),
 \]
for all $n<\omega$ and all $x_0$, \dots, $x_{n-1}$, $y_0$, \dots,
$y_{n-1}\in P$.

\section{Homomorphisms of partial lattices}\label{S:HomsPL}

\begin{definition}\label{D:HomPL}
If $P$ and $Q$ are partial prelattices, a \emph{homomorphism of partial
prelattices} from $P$ to $Q$ is an order-preserving map
$f\colon P\to Q$ such that $a=\bigvee X$
(resp., $a=\bigwedge X$) implies that $f(a)=\bigvee f[X]$ (resp.,
$f(a)=\bigwedge f[X]$), for all $a\in P$ and all $X\in\fine P$. We say
that a homomorphism $f$ is an \emph{embedding}, if $f(a)\leq f(b)$ implies
that $a\leq b$, for all $a$, $b\in P$.
\end{definition}

For a homomorphism $f\colon P\to Q$ of partial lattices, the
\emph{kernel} of $f$, denoted by $\ker f$, is defined as
 \[
 \ker f=\setm{\seq{x,y}\in P\times P}{f(x)\leq f(y)}.
 \]
Moreover, we can define the following maps:

\begin{itemize}
\item The map $\Con f\colon\Con P\to\Con Q$, obtained by defining, for any
congruence $\xa$ of $P$, the congruence $(\Con f)(\xa)$ as the least
congruence of $Q$ that contains all the pairs $\seq{f(x),f(y)}$, for
$\seq{x,y}\in\xa$.

\item The restriction $\Conc f$ of the map $\Con f$ from $\Conc P$ to
$\Conc Q$.

\item The map $\Res f\colon\Con Q\to\Con P$, obtained by defining, for any
congruence $\xb$ of $Q$, the congruence $(\Res f)(\xb)$ as the set of
all $\seq{x,y}\in P\times P$ such that $\seq{f(x),f(y)}\in\xb$. If, in
particular, $P$ is a partial sublattice of $Q$ and $f\colon P\into Q$ is
the inclusion map, then we shall write $\xb\res_P$ instead of
$(\Res f)(\xb)$.

\end{itemize}

This way the maps $P\mapsto\Con P$ and $P\mapsto\Conc P$ can be extended to
functors from partial lattices and their homomorphisms to, respectively,
complete lattices with compactness preserving join-complete
homomorphisms, and \jzs s with \jzh s. On the other hand,
$f\mapsto\Res f$ defines a contravariant functor from partial lattices to
complete lattices with meet-complete homomorphisms that preserve
nonempty directed joins.

The following lemma is a special case of a universal algebraic
triviality:

\begin{lemma}\label{L:CEP}
Let $f\colon P\to Q$ be a homomorphism of partial lattices. Then the
following are equivalent:
\begin{enumerate}
\item $\Con f$ is one-to-one.

\item $\Conc f$ is one-to-one.

\item $\xa=(\Res f)\circ(\Con f)(\xa)$, for all $\xa\in\Con P$.
\end{enumerate}
\end{lemma}

If one of the items of Lemma~\ref{L:CEP} is satisfied, we say that $f$
has the \emph{\cep}.

For a partial lattice $P$, we denote, as in \cite{Wehr}, by $\FL(P)$ the
\emph{free lattice over~$P$}, see \cite{Dean64}. We denote by $j_P$
the canonical embedding from $P$ into $\FL(P)$.

\begin{proposition}\label{P:PtoFLP}
Let $P$ be any partial lattice. Then $j_P$ has the \cep.
\end{proposition}

\begin{proof}
For a congruence $\xa$ of $P$, we denote by $p_\xa$ the canonical
projection from $P$ onto $P/{\xa}$. Since $k=j_{P/{\xa}}\circ p_\xa$
is a homomorphism of partial lattices from $P$ to $\FL(P/{\xa})$,
there exists, by the universal property of the map $j_P$, a unique
lattice homomorphism $q_\xa\colon\FL(P)\onto\FL(P/{\xa})$ such that
$q_\xa\circ j_P=k$, as on the following commutative diagram:
 \[
{
\def\labelstyle{\displaystyle}
\xymatrix{
P\ar@{^(->}[r]^{j_P}\ar@{->>}[d]_{p_\xa}\ar[rd]|-{k} &
\FL(P)\ar@{->>}[d]^{q_\xa} \\
P/{\xa}\ar@{^(->}[r]_{j_{P/{\xa}}} & \FL(P/{\xa})
}
}
 \]
Put $\xb=(\Res j_P)\circ(\Con j_P)(\xa)$, and let $x$, $y\in P$ such
that $x\leq_\xb y$. This means that
$j_P(x)\leq_{(\Con j_P)(\xa)}j_P(y)$, hence, by composing with
$q_\xa$, we obtain that
 \begin{equation}\label{Eq:kxkaky}
 k(x)\leq_{(\Con k)(\xa)}k(y).
 \end{equation}
However, $(\Con k)(\xa)=\Con(j_{P/{\xa}}\circ p_\xa)(\xa)=
(\Con j_{P/{\xa}})(\zero_{P/{\xa}})=\zero_{\FL(P/{\xa})}$. Therefore,
the relation \eqref{Eq:kxkaky} is equivalent to $k(x)\leq k(y)$,
whence, since
$j_{P/{\xa}}$ is an embedding, $p_\xa(x)\leq p_\xa(y)$, that is,
$x\leq_\xa y$. Therefore, $\xb\subseteq\xa$. The converse
inequality is trivial, hence $\xa=\xb$. The conclusion follows.
\end{proof}

Let us recall some further classical definitions, also used in
\cite{Wehr}:

\begin{definition}\label{D:IdFil}
Let $P$ be a partial lattice.
\begin{enumerate}
\item A \emph{partial sublattice} of $P$ is a subset $Q$ of $P$ that is
closed under $\bigvee$ and $\bigwedge$.

\item An \emph{ideal} (resp., \emph{filter}) of $P$ is a lower (resp.,
upper) subset of $P$ closed under $\bigvee$ (resp., $\bigwedge$).
\end{enumerate}
\end{definition}

We observe that both $\es$ and $P$ are simultaneously an ideal and
a filter of~$P$. For a subset $X$ of $P$, we denote by $\II(X)$ (resp.,
$\FF(X)$) the ideal (resp., filter) of $P$ generated by $X$.

\begin{lemma}\label{L:ConfCof}
Let $f\colon P\to Q$ be a homomorphism of partial lattices. If
$\II(f[P])=\FF(f[P])=Q$, then $\Conc f$ is a cofinal map from
$\Conc P$ to $\Conc Q$.
\end{lemma}

\begin{proof}
Put $\xb=(\Con f)(\one_P)$; it suffices to prove that $\xb=\one_Q$.

Fix $x\in P$. Then the relation $f(x)\leq_\xb f(y)$
holds for all $y\in P$, thus the set $F_x=\setm{v\in Q}{f(x)\leq_\xb v}$
contains $f[P]$. Since $F_x$ is obviously a filter of~$Q$, it follows
from the assumptions that $F_x=Q$. Hence, we have established that
 \begin{equation}\label{Eq:f(x)leqv}
 f(x)\leq_\xb v\text{ holds, for all }x\in P\text{ and all }v\in Q.
 \end{equation}
Now it follows from \eqref{Eq:f(x)leqv} that the set
$I_v=\setm{u\in Q}{u\leq_\xb v}$ contains $f[P]$, for all $v\in Q$.
Since $I_v$ is obviously an ideal of $Q$, it follows from the assumptions
that $I_v=Q$. Therefore, $u\leq_\xb v$ for all $u$, $v\in Q$, that is,
$\xb=\one_Q$.
\end{proof}

\begin{corollary}\label{C:FlEmbCof}
Let $P$ be a partial lattice. Then the canonical map\linebreak
$\Conc j_P\colon\Conc P\to\Conc\FL(P)$ is a cofinal embedding.
\end{corollary}

\begin{proof}
By Proposition~\ref{P:PtoFLP}, $\Conc j_P$ is an embedding. Furthermore,
$P$ generates $\FL(P)$ as a lattice, thus, \emph{a fortiori}, $P$
generates $\FL(P)$ both as an ideal and as a filter. Therefore, by
Lemma~\ref{L:ConfCof}, $\Conc j_P$ has cofinal range.
\end{proof}

\section{Duality of complete lattices}\label{S:Dual}

The facts presented in this section are standard, although we do
not know of any reference where they are recorded. Most of the
proofs are straightforward, in which case we omit them. We shall mainly
follow the presentation of \cite{TuWe1}.

In what follows, \emph{\cmh s} are defined in a dual fashion as \cjh s,
and we denote by $\Cj$ (resp., $\Cm$) the category of complete
lattices with \cjh s (resp., \cmh s).

\begin{definition}\label{D:Dual}
Let $A$ and $B$ be complete lattices.
Two maps $f\colon A\to B$ and $g\colon B\to A$ are
\emph{dual}, if the equivalence
 \[
 f(a)\leq b\text{ if and only if }a\leq g(b),
 \]
holds, for all $\seq{a,b}\in A\times B$.
\end{definition}

We recall some basic folklore facts stated in \cite{TuWe1}. For
complete lattices $A$ and $B$, if
$f\colon A\to B$ and $g\colon B\to A$ are dual, then $f$ is a \cjh\
and $g$ is a \cmh. Also, for every \cjh\ (resp., \cmh) $f\colon A\to B$
(resp., $g\colon B\to A$), there exists a unique $g\colon B\to A$
(resp., $f\colon A\to B$) such that $f$ and $g$ are dual, denoted by
$g=f^*$ (resp., $f=g^\dagger$).

The basic categorical properties of the duality thus described
may be recorded in the following lemma.

\begin{lemma}\label{L:CatDual}\hfill
\begin{enumerate}
\item The correspondence $f\mapsto f^*$ defines a contravariant
functor from $\Cj$ to $\Cm$.

\item The correspondence $g\mapsto g^\dagger$ defines a
contravariant functor from $\Cm$ to $\Cj$.

\item If $f$ is a \cjh, then $(f^*)^\dagger=f$.

\item If $g$ is a \cmh, then $(g^\dagger)^*=g$.

\end{enumerate}
\end{lemma}

Of particular importance is the effect of the duality on
\cjh\ of the form $\Con f\colon\Con P\to\Con Q$, where
$f\colon P\to Q$ is a homomorphism of partial lattices.

\begin{lemma}\label{L:DualExt}
Let $P$ and $Q$ be partial lattices, let $f\colon P\to Q$ be a
homomorphism of partial lattices. Then $\Con f$ and $\Res f$ are dual.
\end{lemma}

\begin{lemma}\label{L:gfdagemb}
Let $A$ and $B$ be complete lattices, let $g\colon B\to A$ be a
\cmh. Then $g\circ g^\dagger\circ g=g$. In particular, if $g$ is
surjective, then $g^\dagger$ is an embedding.
\end{lemma}

Let $A$ and $B$ be complete lattices. A map $f\colon A\to B$
is said to \emph{preserve nonempty directed joins}, if
$f\left(\bigvee X\right)=\bigvee f[X]$, for any nonempty upward directed
subset $X$ of $A$.

\begin{lemma}\label{L:CompPres}
Let $A$ and $B$ be complete lattices. Then the following assertions hold:
\begin{enumerate}
\item Let $g\colon B\to A$ be a \cmh. If $g$ preserves nonempty directed
joins, then the dual map $g^\dagger\colon A\to B$ preserves compactness.

\item Let $f\colon A\to B$ be a \cjh. If $A$ is algebraic and $f$
preserves compactness, then the dual map $f^*\colon B\to A$ preserves
nonempty directed joins.
\end{enumerate}
\end{lemma}

\begin{proof}
(i) Let $a\in\KK(A)$, we prove that $b=g^\dagger(a)$ belongs to
$\KK(B)$. So let $X$ be a nonempty upward directed subset of $B$ such that
$b\leq\bigvee X$. By the definition of $g^\dagger$, this means that
$a\leq g\left(\bigvee X\right)$, which, by the assumption on~$g$, can be
written $a\leq\bigvee g[X]$. Therefore, since $a\in\KK(A)$, there exists
$x\in X$ such that $a\leq g(x)$, that is, $b\leq x$. Hence
$b\in\KK(B)$.

(ii) Let $Y$ be an upward directed subset of $B$, put $b=\bigvee Y$. Let
$a\in\KK(A)$ such that $a\leq f^*(b)$. This means that $f(a)\leq b$,
however, $f(a)$ is, by assumption on $f$, compact in $B$, thus
$f(a)\leq y$ for some $y\in Y$, whence $a\leq\bigvee f^*[Y]$. Since $A$
is algebraic, this proves that $f^*(b)\leq\bigvee f^*[Y]$. The converse
inequality is trivial.
\end{proof}

As a corollary, we get the following well-known fact, see, \emph{e.g.},
Lemma~1.3.3 in~\cite{Gorb}:

\begin{lemma}\label{L:ClSAlg}
let $A$ be an algebraic lattice, let $B$ be a \emph{closure system}
in~$A$, \emph{i.e.}, a complete meet-subsemilattice of $A$ that is closed
under nonempty directed joins. Then $B$ is an algebraic lattice.
\end{lemma}

\begin{proof}
Let $g\colon B\into A$ be the inclusion map. By assumption and by
Lem\-ma~\ref{L:CompPres}(i), the dual map $f=g^\dagger$ preserves
compactness. Let $b\in B$. For any $x\in\KK(B)$ such that $x\leq b$,
the inequalities $x\leq f(x)\leq f(b)=b$ hold, whence
 \[
 b=\bigvee\setm{x\in\KK(A)}{x\leq b}=
 \bigvee\setm{f(x)}{x\in\KK(A),\ x\leq b}.
 \]
The conclusion follows from the fact that $f[\KK(A)]\subseteq\KK(B)$.
\end{proof}

\section{\Ccb\ semilattices}\label{S:LCB}

\begin{definition}\label{D:CB}
Let $S$ be a \jzs. We say that $S$ is
\begin{itemize}
\item \emph{\cb}, if $S$ is a complete lattice and it satisfies the
infinite meet distributivity law (MID), that is, it satisfies the
infinitary identity
 \begin{equation}
 a\vee\bigwedge_{i\in I}x_i=\bigwedge_{i\in I}(a\vee x_i),\tag{MID}
 \end{equation}
where $a$ and the $x_i$-s range over the elements of $S$.

\item \emph{\ccb}, if every principal ideal of $S$ is \cb.

\end{itemize}
\end{definition}

Equivalently, $S$ is \cb\ if{f} $S$ is a dually relatively
\pscp\ complete lattice, see \cite{GLT2} for explanation
about the latter terminology.

We observe that every \ccb\ lattice is, of course, distributive.

The crucial point that we shall use about \ccb\ lattices is
the following:

\begin{lemma}\label{L:InjLcb}
Let $S$ be a \ccb\ lattice, let $A$ be a cofinal \jz-subsemilattice of
a \jzs\ $B$. Then every \jzh\ from $A$ to $S$ extends to some \jzh\
from $B$ to~$S$.
\end{lemma}

\begin{proof}
The conclusion of Lemma~\ref{L:InjLcb} follows immediately from
Theorem~3.11 of~\cite{Wehr92}. However, it is worth observing
that since we are dealing with semilattices, there is also a direct
proof. Namely, if $f\colon A\to S$ is any \jzh, the
completeness assumption on $S$ and the fact that $f$ has cofinal range
make it possible to define a map $g\colon B\to S$ by the rule
 \[
 g(b)=\bigwedge\setm{f(x)}{x\in A\text{ and }b\leq x}.
 \]
It follows then from (MID) that $g$ is a join homomorphism. It is
obvious that $g$ extends $f$.
\end{proof}

\begin{remark}
By using some of the techniques of the proof of Theorem~3.11 of
\cite{Wehr92}, it is not hard to prove that in fact,
Lemma~\ref{L:InjLcb} \emph{characterizes} \ccb\ lattices.
\end{remark}

Now we can already provide a proof of Theorem~\ref{T:1Main} stated in the
Introduction.

\begin{proof}[Proof of Theorem~\textup{\ref{T:1Main}}]
By Corollary~\ref{C:FlEmbCof} and Lemma~\ref{L:InjLcb}, there exists a
\jzh\ $\psi\colon\Conc\FL(P)\to S$ such that
$\psi\circ\Conc j_P=\varphi$. Then it suffices to apply
Theorem~\ref{T:1DimLift} to $\psi$.
\end{proof}

Now we can provide a proof of Theorem~\ref{T:2Main} stated in the
Introduction:

\begin{proof}[Proof of Theorem~\textup{\ref{T:2Main}}]
Let $\DD$ be described by homomorphisms $f\colon K\to P$ and $g\colon K\to
Q$ of partial lattices, with $K$ a lattice, and let $\varphi$ be described
by \jzh s $\mu\colon\Conc P\to S$ and $\nu\colon\Conc Q\to S$ such that
$\mu\circ\Conc f=\nu\circ\Conc g$. We shall construct a relatively
complemented lattice $L$, homomorphisms $\oll{f}\colon P\to L$ and
$\ol{g}\colon Q\to L$ of partial lattices, and an isomorphism
$\eps\colon\Conc L\to S$ such that $\oll{f}\circ f=\ol{g}\circ g$,
$\mu=\eps\circ\Conc\oll{f}$, $\nu=\eps\circ\Conc\ol{g}$, $L$ is
relatively complemented, $\oll{f}[P]\cup\ol{g}[Q]$ generates $L$ as an
ideal (resp., filter), and, if $S$ is generated as an ideal by
$\rng\mu\cup\rng\nu$, then $L$ is generated by
$\oll{f}[P]\cup\ol{g}[Q]$ as a convex sublattice (where $\rng\mu$ stands
for the range of $\mu$).

We first reduce the problem to the case where both $f$ and $g$ are
\emph{embeddings}, as follows (see also the end of the proof of
Proposition~18.5 of \cite{Wehr}). We put
$\lambda=\mu\circ\Conc f=\nu\circ\Conc g$, and we define congruences
$\xd\in\Con K$, $\xa\in\Con P$, and $\xb\in\Con Q$ as follows:
 \begin{align*}
 \xd&=\setm{\seq{x,y}\in K\times K}{\lambda\Theta_K^+(x,y)=0},\\
 \xa&=\setm{\seq{x,y}\in P\times P}{\mu\Theta_P^+(x,y)=0},\\
 \xb&=\setm{\seq{x,y}\in Q\times Q}{\nu\Theta_Q^+(x,y)=0}.
 \end{align*}
We denote by $p_\xd\colon K\onto K/{\xd}$,
$p_\xa\colon P\onto P/{\xa}$, $p_\xb\colon Q\onto Q/{\xb}$ the
canonical projections. Then there are unique homomorphisms of partial
lattices $f'\colon K/{\xd}\into P/{\xa}$ and
$g'\colon K/{\xd}\into Q/{\xb}$ such that
$f'\circ p_{\xd}=p_{\xa}\circ f$ and
$g'\circ p_{\xd}=p_{\xb}\circ g$, and both $f'$ and $g'$ are
embeddings. Furthermore, we can define \jzh s
$\mu'\colon\Conc(P/{\xa})\to S$ and $\nu'\colon\Conc(Q/{\xb})\to S$
by the rules $\mu'(\xx\vee\xa/{\xa})=\mu(\xx)$ for all
$\xx\in\Conc P$, and $\nu'(\xx\vee\xb/{\xb})=\nu(\xx)$ for all
$\xx\in\Conc Q$.

Since $\mu'\circ\Conc f'=\nu'\circ\Conc g'$ and both
$f'$ and $g'$ are embeddings, there are, by assumption, a relatively
complemented lattice $L$, homomorphisms
$\oll{f'}\colon P/{\xa}\to L$ and
$\ol{g'}\colon Q/{\xb}\to L$ of partial lattices, and an isomorphism
$\eps\colon\Conc L\to S$ such that $\oll{f'}\circ f'=\ol{g'}\circ g'$,
$\mu'=\eps\circ\Conc\oll{f'}$, $\nu'=\eps\circ\Conc\ol{g'}$, $L$ is
relatively complemented, $\oll{f'}[P/{\xa}]\cup\ol{g'}[Q/{\xb}]$
generates $L$ as an ideal (resp., filter), and, if $S$ is generated as
an ideal by $\rng\mu'\cup\rng\nu'$, then $L$ is generated by
$\oll{f'}[P/{\xa}]\cup\ol{g'}[Q/{\xb}]$ as a convex sublattice. Then
$f'=\oll{f'}\circ p_\xa$ and $g'=\ol{g'}\circ p_\xb$, together with
$\eps$ and $L$, solve the amalgamation problem for $f$ and $g$.

Hence we can reduce the problem to the case where
both $f$ and $g$ are embeddings. Without loss of generality, $f$
and $g$ are, respectively, the set-theoretical inclusion from $K$ into
$P$ (resp., $Q$), and $K=P\cap Q$.

Then we define a partial lattice $R$ as follows, whose classical
construction is also recalled in the statement of Proposition~3.4 in
\cite{Wehr}. The underlying set of $R$ is $P\cup Q$, and the partial
ordering of $R$ is defined as follows. For $x$, $y\in R$, the inequality
$x\leq y$ holds if{f} one of the following cases hold:
\begin{enumerate}
\item $x$, $y\in P$ and $x\leq y$ in $P$;

\item $x$, $y\in Q$ and $x\leq y$ in $Q$;

\item $x\in P$, $y\in Q$, and there exists $z\in K$ such that $x\leq z$
in $P$ and $z\leq y$ in $Q$.

\item $x\in Q$, $y\in P$, and there exists $z\in K$ such that $x\leq z$
in $Q$ and $z\leq y$ in $P$.

\end{enumerate}

The partially ordered set $R$ can be given a structure of partial
lattice, as follows. For $a\in R$ and $X\in\fine{R}$, $a=\bigvee X$
holds in $R$, if either $X\cup\{a\}\subseteq P$ and $a=\bigvee X$ in $P$
or $X\cup\{a\}\subseteq Q$ and $a=\bigvee X$ in $Q$. The meet operation
on $R$ is defined dually.

Let $u$ (resp., $v$) be the inclusion map from $P$ (resp., $Q$) into
$R$. It is stated in \cite{Wehr}, and very easy to prove, that
$\seq{R,u,v}$ is a \emph{pushout} of the diagram $\seq{K,P,Q,f,g}$ in the
category of partial lattices and their homomorphisms. We shall abuse the
notation by stating this as $R=P\amalg_KQ$, the maps $f$ and $g$ then
being understood.

Now we put
$\xC=\setm{\seq{\xa,\xb}\in\Con P\times\Con Q}{\xa\res_K=\xb\res_K}$. It
is obvious that $\xC$ is a complete meet-subsemilattice of
$\Con P\times\Con Q$, closed under nonempty directed
suprema. Hence, by Lemma~\ref{L:ClSAlg}, $\xC$ is an algebraic lattice.
Observe that $\seq{\zero_P,\zero_Q}\in\xC$.

Let $\varphi\colon\Con R\to\xC$, $\xc\mapsto\seq{\xc\res_P,\xc\res_Q}$.
Then $\varphi$ is a \cmh, and it preserves nonempty directed joins.
Hence, by Lemma~\ref{L:CompPres}(i), the dual map $\psi=\varphi^\dagger$
of $\varphi$ is a compactness-preserving \cjh\ from $\xC$ to $\Con R$.

\setcounter{claim}{0}
\begin{claim}\label{Cl:psiEmb}
The map $\varphi$ is surjective, while $\psi$ is an embedding.
\end{claim}

\begin{cproof}
By Lemma~\ref{L:gfdagemb}, it suffices to prove that $\varphi$ is
surjective. Let $\seq{\xa,\xb}\in\xC$, put $\xd=\xa\res_K=\xb\res_K$.
Then the natural homomorphism from $K/\xd$ into $P/{\xa}$ (resp.,
$Q/{\xb}$) is an embedding, therefore, by using the universal property
of $R=P\amalg_KQ$, there exists a homomorphism
$r\colon P\amalg_KQ\onto(P/{\xa})\amalg_{(K/\xd)}(Q/{\xb})$ such that
the following diagram commutes ($p_\xa$ and $q_\xb$ denote the
canonical projections):
 \[
{
\def\labelstyle{\displaystyle}
\xymatrixcolsep{1pc}
\xymatrix{
 & (P/{\xa})\amalg_{(K/\xd)}(Q/{\xb}) & \\
P/{\xa}\ar@{_(->}[ru]^{} & P\amalg_KQ\ar@{->>}[u]^{r} &
Q/{\xb}\ar@{^(->}[lu]^{} \\
P\ar@{->>}[u]^{p_\xa}\ar@{_(->}[ru]^(.3){u} &
K/\xd\ar@{_(->}[lu]^{}\ar@{^(->}[ru]^{} &
Q\ar@{->>}[u]_{q_\xb}\ar@{^(->}[lu]_(.3){v}\\
& K\ar@{_(->}[lu]^{f}\ar@{->>}[u]^{}\ar@{^(->}[ru]_{g} &
}
}
 \]
Put $\xc=\ker r$. Then $\xc$ is a congruence of $R$. Moreover, for any
$x$, $y\in P$, $x\leq_\xc y$ if{f} $r(x)\leq r(y)$, that is,
$p_\xa(x)\leq p_\xa(y)$, or $x\leq_\xa y$. Hence $\xc\res_P=\xa$.
Similarly, $\xc\res_Q=\xb$, hence $\varphi(\xc)=\seq{\xa,\xb}$.
\end{cproof}

\begin{claim}\label{Cl:K(f)}
The map $\KK(\psi)$ is cofinal from $\KK(\xC)$ to $\Conc R$.
\end{claim}

\begin{cproof}
Put $\xc=\psi(\seq{\one_P,\one_Q})$. It follows from the definition of
$\psi$ that $\xc\res_P=\one_P$ and $\xc\res_Q=\one_Q$. Pick $z\in K$
(we have supposed that $K\neq\es$). For any $x\in P$ and $y\in Q$,
$x\leq_\xc z$ (because $\xc\res_P=\one_P$) and $z\leq_\xc y$ (because
$\xc\res_Q=\one_Q$), hence $x\leq_\xc y$. Similarly, $y\leq_\xc x$.
Therefore, $\psi(\seq{\one_P,\one_Q})=\xc=\one_R$. The conclusion of
Claim~\ref{Cl:K(f)} follows.
\end{cproof}

\begin{claim}\label{Cl:Apush}
$\KK(\xC)$ is a pushout of $\Conc P$ and $\Conc Q$ above $\Conc f$ and
$\Conc g$ in the category of all \jz-semilattices.
\end{claim}

\begin{cproof}
Let $\xi\colon\xC\to\Con P$ and $\eta\colon\xC\to\Con Q$ be the
canonical projections. By the definition of $\xC$, the following diagram
 \[
{
\def\labelstyle{\displaystyle}
\xymatrixcolsep{1pc}
\xymatrix{
 & \xC\ar[dl]_{\xi}\ar[dr]^{\eta} & \\
\Con P\ar[dr]_{\Res f} & & \Con Q\ar[dl]^{\Res g}\\
 & \Con K &
}
}
 \]
is a \emph{pullback} in the category of all algebraic lattices with \cmh s
that preserve nonempty directed joins. By dualizing this diagram (see
Lemmas \ref{L:CatDual}, \ref{L:DualExt}, and \ref{L:CompPres}), then
by taking the image of the new diagram under the functor $\KK$, and then
by using Proposition~\ref{P:AlgLattSem}, we obtain successively the
following diagrams, the left hand side a pushout in the category of all
algebraic lattices and compactness preserving \cjh s, the right hand
side a pushout in the category of \jzs s with \jzh s,
 \[
{
\def\labelstyle{\displaystyle}
\xymatrixcolsep{1pc}
\xymatrix{
 & \xC & & & \KK(\xC) &\\
 \Con P\ar[ur]^{\xi^\dagger} & &  \Con Q\ar[ul]_{\eta^\dagger} &
 \Conc P\ar[ur]^{\alpha=\KK(\xi^\dagger)}
 & & \Conc Q\ar[ul]_{\beta=\KK(\eta^\dagger)} \\
 & \Con K\ar[ul]^{\Con f}\ar[ur]_{\Con g} & & &
 \Conc K\ar[ul]^{\Conc f}\ar[ur]_{\Conc g} &
}
}
 \]
which completes the proof of Claim~\ref{Cl:Apush}.
\end{cproof}

By applying the sequence of two functors used in the proof of
Claim~\ref{Cl:Apush} to the following commutative diagram,
 \[
{
\def\labelstyle{\displaystyle}
\xymatrixcolsep{1pc}
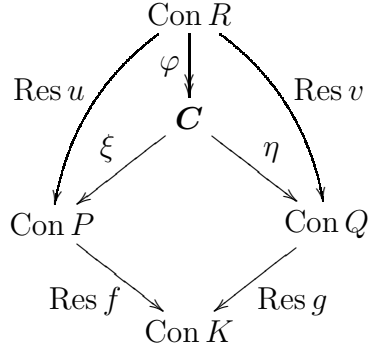
\xymatrix{
 & \Con R\ar@/_1pc/[ddl]_{\Res u}\ar@/^1pc/[ddr]^{\Res v}
 \ar@{->>}[d]_{\varphi} & \\
 & \xC\ar[dl]_{\xi}\ar[dr]^{\eta} & \\
\Con P\ar[dr]_{\Res f} & & \Con Q\ar[dl]^{\Res g}\\
 & \Con K &
}
}
 \]
we obtain, successively, the two following commutative diagrams:
 \[
{
\def\labelstyle{\displaystyle}
\xymatrixcolsep{1pc}
\xymatrix{
 & \Con R & & & \Conc R & \\
 & \xC\ar@{_(->}[u]^{\psi} & & &
 \KK(\xC)\ar@{_(->}[u]|-{\strut\KK(\psi)} &\\
 \Con P\ar[ur]^{\xi^\dagger}\ar@/^1pc/[uur]^{\Con u}
 & &  \Con Q\ar[ul]_{\eta^\dagger}\ar@/_1pc/[uul]_{\Con v} &
 \Conc P\ar[ur]^{\alpha}\ar@/^1pc/[uur]^{\Conc u}
 & & \Conc Q\ar[ul]_{\beta}\ar@/_1pc/[uul]_{\Conc v} \\
 & \Con K\ar[ul]^{\Con f}\ar[ur]_{\Con g} & & &
 \Conc K\ar[ul]^{\Conc f}\ar[ur]_{\Conc g} &
}
}
 \]
Since $\mu\circ\Conc f=\nu\circ\Conc g$ and by Claim~\ref{Cl:Apush},
there exists a \jzh\ $\gamma\colon\KK(\xC)\to S$ such that the following
diagram
 \[
{
\def\labelstyle{\displaystyle}
\xymatrixcolsep{1pc}
\xymatrix{
 & S & \\
 & \KK(\xC)\ar[u]^{\gamma} & \\
 \Conc P\ar[ur]^{\alpha}\ar@/^1pc/[uur]^{\mu}
 & & \Conc Q\ar[ul]_{\beta}\ar@/_1pc/[uul]_{\nu} \\
 & \Conc K\ar[ul]^{\Conc f}\ar[ur]_{\Conc g} &
}
}
 \]
is commutative. Furthermore, by Claims \ref{Cl:psiEmb} and
\ref{Cl:K(f)}, $\KK(\psi)$ is a cofinal embedding from $\KK(\xC)$ into
$\Conc R$, while, by Corollary~\ref{C:FlEmbCof}, $\Conc j_R$ is a cofinal
embedding from $\Conc R$ into $\Conc\FL(R)$. Therefore, the map
$(\Conc j_R)\circ\KK(\psi)$ is a cofinal embedding from $\KK(\xC)$ into
$\Conc\FL(R)$. By Lemma~\ref{L:InjLcb}, there exists a \jzh\
$\pi\colon\Conc\FL(R)\to S$ such that\linebreak
$\pi\circ(\Conc j_R)\circ\KK(\psi)=\gamma$. By Theorem~\ref{T:1DimLift},
there are a relatively complemented lattice $L$, a lattice homomorphism
$h\colon\FL(R)\to L$, and an isomorphism $\eps\colon\Conc L\to S$ such
that $\pi=\eps\circ\Conc h$, the range of $h$ generates $L$ as an ideal
(resp., filter), and, if the range of $\pi$ is cofinal in
$S$, then the range of $h$ generates $L$ as a convex sublattice. The
latter condition is certainly satisfied if
$\rng\mu\cup\rng\nu$ is cofinal in~$S$ (because $\rng\gamma$ contains
$\rng\mu\cup\rng\nu$). Some of this information is summarized on the
following commutative diagram.
 \[
{
\def\labelstyle{\displaystyle}
\xymatrixcolsep{1pc}
\xymatrix{
 S \\
 \KK(\xC)\ar@{^(->}[r]_{\KK(\psi)}\ar[u]^{\gamma} &
 \Conc R\ar@{^(->}[0,2]_{\Conc j_R} & &
 \Conc\FL(R)\ar[0,2]_{\Conc h}\ar[-1,-3]^(.6){\pi} & &
 \Conc L\ar@/_1pc/[-1,-5]_{\eps}
}
}
 \]
Now we consider the following commutative diagram:
 \[
{
\def\labelstyle{\displaystyle}
\xymatrixcolsep{1pc}
\xymatrix{
 & \FL(R) & \\
 P\ar[ur]^{f'=j_R\circ u} & & Q\ar[ul]_{g'=j_R\circ v}\\
 & K\ar[ul]^f\ar[ur]_g &
}
}
 \]
We further compute:
 \[
 \pi\circ\Conc f'=\pi\circ\Conc j_R\circ\Conc u=
 \pi\circ\Conc j_R\circ\KK(\psi)\circ\alpha
 =\gamma\circ\alpha=\mu.
 \]
A similar argument proves the equality
$\pi\circ\Conc g'=\nu$. The fact that $\FL(R)$ is generated, as
a lattice, by $f'[P]\cup g'[Q]$, trivially follows from $R=P\cup Q$.
Therefore, the maps $\oll{f}=h\circ f'$ and $\ol{g}=h\circ g'$, together
with the isomorphism~$\eps$, satisfy the required conditions.
\end{proof}

The following corollary generalizes Theorem~2 of \cite{GLW}:

\begin{corollary}\label{C:MoreRepr}
Let $S$ be a distributive \jzs\ that can be expressed as a \jz-direct
limit of at most $\aleph_1$ \ccb\ lattices. Then there exists a
relatively complemented lattice $L$ with zero such that $\Conc L\cong S$.
Furthermore, if $S$ has a largest element, then $L$ can be taken
bounded.
\end{corollary}

\begin{proof}
Write $S=\varinjlim(S_i)_{i\in I}$ with transition \jzh s\linebreak
$f_{i,j}\colon S_i\to S_j$, limiting maps $f_i\colon S_i\to S$, where
$I$ is an upward directed partially ordered set of size at most $\aleph_1$
and all the $S_i$-s are \ccb\ lattices. As at the beginning of the proof of
Theorem~2 of \cite{GLW}, we may assume without loss of generality that $I$
is a \emph{$2$-ladder}, that is, a lattice with zero in which every
principal ideal is finite and every element has at most
two immediate predecessors. The rest of the proof goes as
the proof of Theorem~2 of
\cite{GLW}, by using Theorem~\ref{T:2Main} for the amalgamation step.
\end{proof}

The following corollary generalizes Theorem~3 of \cite{GLW}. Its proof
is similar, again by using Theorem~\ref{T:2Main}.

\begin{corollary}\label{C:CPE}
Let $K$ be a lattice that can be expressed as a direct union of
\emph{countably} many lattices whose congruence semilattices are \ccb.
Then $K$ embeds congruence-preservingly into some relatively
complemented lattice~$L$, which it generates as a convex sublattice.
\end{corollary}

\section{Conditionally $\kappa$-co-Brouwerian semilattices}\label{S:LockBr}

\begin{definition}\label{D:kCB}
Let $S$ be a \jzs, let $\kappa$ be an infinite cardinal. We say that
$S$ is \emph{\ckcb}, if it satisfies the following conditions.
\begin{enumerate}
\item \emph{$<\kappa$-interpolation property}: for all nonempty $X$,
$Y\subseteq S$ such that $|X|,|Y|<\kappa$ such that $X\leq Y$ (that
is, $x\leq y$ for all $\seq{x,y}\in X\times Y$), there exists $z\in S$
such that $X\leq z\leq Y$.

\item \emph{$<\kappa$-interval axiom}: for all $X\subseteq S$ such that
$|X|<\kappa$ and all $a$, $b\in S$ such that $a\leq b\vee x$ for all
$x\in X$, there exists $c\in S$ such that $a\leq b\vee c$ and $c\leq X$.
\end{enumerate}
\end{definition}

Observe that every \ckcb\ \jzs\ is obviously distributive (take $X$ a
pair in (ii)).

Now we prove the following analogue of Lemma~\ref{L:InjLcb}:

\begin{lemma}\label{L:InjLkcb}
Let $\kappa$ be an infinite cardinal, let $S$ be a \ckcb\ \jzs, let
$A$ be a cofinal \jz-subsemilattice of a \jzs\ $B$ such that
$|B|<\kappa$. Then every \jzh\ from $A$ to $S$ extends to some \jzh\
from $B$ to $S$.
\end{lemma}

\begin{proof}
It suffices to consider the case where $B$ is a monogenic extension of
$A$, \emph{i.e.}, $B=A[b]=A\cup\setm{x\vee b}{x\in A}$, where $b$ is an
element of $B$. Let $f\colon A\to S$ be a \jzh. Let
$\setm{\seq{x_i,y_i}}{i\in I}$ enumerate all elements $\seq{x,y}$ of
$A\times A$ such that $x\leq y\vee b$, and let $\setm{z_j}{j\in J}$
enumerate all elements $z$ of $A$ such that $b\leq z$, with
$|I|,|J|<\kappa$. Observe that $I\neq\es$, and, since $A$ is cofinal in
$B$, $J\neq\es$. For all $\seq{i,j}\in I\times J$, the inequality
$x_i\leq y_i\vee z_j$ holds, thus $f(x_i)\leq f(y_i)\vee f(z_j)$. By the
$<\kappa$-interval axiom, for all $i\in I$, there exists $\xb_i\in S$
such that $f(x_i)\leq f(y_i)\vee\xb_i$ and $\xb_i\leq f(z_j)$ for all
$j\in J$. By the $<\kappa$-interpolation property, there exists
$\xb\in S$ such that $\xb_i\leq\xb\leq f(z_j)$ for all
$\seq{i,j}\in I\times J$. Hence, $f(x_i)\leq f(y_i)\vee\xb$ for all
$i\in I$, and $\xb\leq f(z_j)$ for all $j\in J$, so that there exists a
unique \jzh\ $g\colon B\to S$ extending $f$ such that $g(b)=\xb$.
\end{proof}

We shall now outline a proof of the following analogue of
Proposition~18.5 of~\cite{Wehr}.

\begin{lemma}\label{L:Factkappa}
Let $\kappa$ be an uncountable cardinal, let $\DD$ be a truncated square
of partial lattices, with lattice bottom, of size less than $\kappa$,
let $S$ be a \ckcb\ \jzs. Then every homomorphism
$\varphi\colon\Conc\DD\to S$ has a factor of the form
$f\colon\DD\to L$, where $L$ is a lattice generated by the range of
$f$ (thus, $|L|<\kappa$).
\end{lemma}

\begin{proof}
We use the same notation as for the proof of Theorem~\ref{T:2Main}
presented in Section~\ref{S:LCB}. In particular, $|K|,|P|,|Q|<\kappa$.
Then the proof of Theorem~\ref{T:2Main} applies \emph{mutatis mutandis},
by using Lemma~\ref{L:InjLkcb} instead of Lemma~\ref{L:InjLcb}, to
establish that the canonical pushout homomorphism $f\colon\DD\to\FL(R)$,
with $R=P\amalg_KQ$, is a factor of $\varphi$: all semilattices that need
to be of size less than $\kappa$ are indeed of size less than $\kappa$,
moreover, the last extension step from $\FL(R)$ to $L$ used in the proof
of Theorem~\ref{T:2Main} is no longer necessary since we require only
`factor' instead of `lift'. Observe that since $\kappa$ is uncountable,
$L=\FL(R)$ still has size less than $\kappa$.
\end{proof}

Our next definitions are borrowed from \cite{Wehr}:

\begin{definition}\label{D:D0PL}
Let $S$ be a \jzs. A \emph{\MPL{S}} is a pair $\seq{P,\mu}$,
where $P$ is a partial lattice and $\mu\colon\Conc P\to S$ is a
\jzh. If, in addition, $P$ is a lattice, we say that $\seq{P,\mu}$ is a
\emph{\ML{S}}.

A \MPL{S} $\seq{P,\mu}$ is \emph{proper}, if $\mu$ \emph{isolates
zero}, that is, $\mu^{-1}\{0\}=\{\zero_P\}$.
\end{definition}

\begin{definition}\label{D:EmbDPL}
Let $S$ be a \jzs, let $\seq{P,\mu}$ and
$\seq{Q,\nu}$ be \MPL{S}s. A \emph{homomorphism} from $\seq{P,\mu}$ to
$\seq{Q,\nu}$ is a homomorphism $f\colon P\to Q$ of partial
lattices such that $\nu\circ\Conc f=\mu$. If, in addition, $f$ is an
embedding of partial lattices, we say that $f$ is an \emph{embedding} of
\MPL{S}s.
\end{definition}

\begin{definition}\label{D:ExtMPL}
Let $S$ be a \jzs,
let $\seq{P,\mu}$ and $\seq{L,\varphi}$ be \MPL{S}s, with $L$ a lattice.
We say that an embedding $f\colon\seq{P,\mu}\into\seq{L,\varphi}$ is
a \emph{lower embedding} (resp., \emph{upper embedding}, \emph{internal
embedding}), if the filter (resp., ideal, convex sublattice) of $L$
generated by $P$ equals~$L$.
\end{definition}

\begin{definition}\label{D:EqnSat}
Let $S$ be a \jzs, let $X$ be a subset of $S$.
A proper \ML{S} $\seq{L,\varphi}$ is \emph{$X$-saturated}
(resp., \emph{lower $X$-saturated}, \emph{upper $X$-saturated},
\emph{internally $X$-saturated}), if for every embedding (resp., lower
embedding, upper embedding, internal embedding)
$e\colon\seq{K,\lambda}\into\seq{P,\mu}$ of finite proper \MPL{S}s
such that $\rng\mu\subseteq X\cup\rng\varphi$,
with $K$ a \emph{lattice}, and every homomorphism
$f\colon\seq{K,\lambda}\to\seq{L,\varphi}$, there exists a homomorphism
$g\colon\seq{P,\mu}\to\seq{L,\varphi}$ such that $g\circ e=f$.
\end{definition}

Now a standard increasing chain argument makes it possible to prove the
following result.

\begin{proposition}\label{P:EmbSat}
Let $\kappa$ be an uncountable cardinal, let $S$ be a \ckcb\ \jzs, let
$X\subseteq S$ such that $|X|<\kappa$. Every proper \MPL{S}
$\seq{P,\varphi}$ such that $|P|<\kappa$ admits an embedding (resp., a
lower embedding, an upper embedding, an internal embedding) into a
$X$-saturated (resp.,  lower $X$-saturated, upper $X$-saturated,
internally $X$-saturated) \ML{S} $\seq{L,\psi}$ such that
$|L|=|P|+|X|+\aleph_0$.
\end{proposition}

\begin{proof}
As in the proof of Proposition 19.3 of \cite{Wehr}. We first use
Corollary~\ref{C:FlEmbCof} and Lemma~\ref{L:InjLkcb} to extend
$\seq{P,\varphi}$ by $\seq{\FL(P),\psi}$ for some $\psi$. Then the
\MPL{S} $\seq{\FL(P),\psi}$ may not be proper, so we need to replace
it by its quotient under the congruence of $\FL(P)$ that consists of
all pairs $\seq{x,y}$ such that $\psi\Theta^+(x,y)=0$ (called the
\emph{kernel projection} in \cite{Wehr}).

This way, we obtain that $P$ may be assumed to be a lattice from the
start. Furthermore, there are at most $|P|+|X|+\aleph_0$ pairs of the
form $\seq{e,f}$ where $e\colon\seq{K,\lambda}\to\seq{Q,\nu}$ and
$f\colon\seq{K,\lambda}\to\seq{P,\varphi}$ are homomorphisms of
\MPL{S}s with $K$ a lattice, both $K$ and $Q$ finite, $e$ an
embedding, and $\rng\nu$ is contained in $X\cup\rng\varphi$. We
increase $\seq{P,\varphi}$ by a transfinite sequence of length
$|P|+|X|+\aleph_0$ of \ML{S}s. At each stage $\seq{L,\psi}$ of the
construction, we pick the corresponding pair $\seq{e,f}$ of
homomorphisms. The amalgamation result of Lemma~\ref{L:Factkappa}
makes it possible to find a \ML{S} $\seq{L',\psi'}$, together with
homomorphisms $e'$ and $f'$, such that the following diagram commutes:
 \[
 {
 \def\labelstyle{\displaystyle}
 \xymatrix{
 \seq{Q,\nu}\ar[r]^{f'} & \seq{L',\psi'}\\
 \seq{K,\lambda}\ar[u]^{e}\ar[r]_{f} &
 \seq{L,\psi}\ar[u]_{e'}
 }
 }
 \]
Again, by replacing $\seq{L',\psi'}$ by its quotient under its kernel
projection, we may assume that $\seq{L',\psi'}$ is proper; let, then,
$\seq{L',\psi'}$ be the next step of the construction.

We denote by $\seq{P,\varphi}^*$ the direct limit of that
construction. Iterating $\omega$ times the operation
$\seq{P,\varphi}\into\seq{P,\varphi}^*$ and taking again the direct
limit yields the desired result.
\end{proof}

Now, by using Proposition~\ref{P:EmbSat}, we argue as in Section~20 of
\cite{Wehr} to obtain the following analogue of Proposition~20.8 of
\cite{Wehr}. Observe that the proof is, in fact, much simpler than the
one of Proposition~20.8 of \cite{Wehr}. The reason for this is that we no
longer need to check that the corresponding \MPL{S}s are `balanced',
which removes lots of technical complexity.

\begin{proposition}\label{P:SummSat}
Let $\kappa$ be an uncountable cardinal, let $S$ be a \ckcb\ \jzs,
let $\seq{L,\varphi}$ be an internally $X$-saturated \MPL{S}. Then the
following assertions hold:
\begin{enumerate}
\item $L$ is relatively complemented.

\item The map $\varphi$ is an embedding from $\Conc L$ into $S$, and
$X\cap\dnw\rng\varphi\subseteq\rng\varphi$. (For a subset $Y$ of $S$,
$\dnw Y$ denotes the lower subset of $S$ generated by $Y$.)

\item For $o$, $a$, $b$, $i\in L$ such that $o\leq\set{a,b}\leq i$,
$\Theta_L(o,a)=\Theta_L(o,b)$ if{f} there are $a_0$, $a_1$, $b_0$,
$b_1\in[o,i]$ such that the following conditions hold:
\begin{enumerate}
\item $a=a_0\vee a_1$, $b=b_0\vee b_1$, and
$a_0\wedge a_1=b_0\wedge b_1=o$.

\item $a_0$ and $b_0$ (resp., $a_1$ and $b_1$) are perspective in
$[o,i]$, \emph{i.e.}, for all $l<2$, there exists $x\in[o,i]$ such that
$a_l\wedge x=b_l\wedge x=o$ and $a_l\vee x=b_l\vee x=i$.
\end{enumerate}

\item If, in addition, $\seq{L,\varphi}$ is either lower $X$-saturated or
upper $X$-sat\-u\-rat\-ed, then $X\subseteq\rng\varphi$.
\end{enumerate}
\end{proposition}

\begin{proof}[Outline of Proof]
We imitate the proof of Proposition~20.8 of \cite{Wehr}. We first
show, for example, that $L$ is relatively complemented. For $a<b<c$
in $L$, let $K=\set{a,b,c}$ be the three-element chain, let
$f\colon K\into L$ be the natural embedding, and put
$\lambda=\varphi\circ\Conc f$. Then $\seq{K,\lambda}$ is a
finite, proper \ML{S} and $f$ is an embedding from $\seq{K,\lambda}$
into $\seq{L,\varphi}$. Next, we put $P=\set{a,b,c,t}$, the Boolean
lattice with bottom $a$, top $c$, and atoms $b$ and $t$, endowed with
the homomorphism $\mu\colon\Conc P\to S$ defined by
 \begin{align*}
 \mu\Theta_P(a,b)=\mu\Theta_P(t,c)&=\varphi\Theta_L(a,b),\\
 \mu\Theta_P(a,t)=\mu\Theta_P(b,c)&=\varphi\Theta_L(b,c).
 \end{align*}
Then $\seq{P,\mu}$ is a proper \ML{S}, with
$\rng\mu\subseteq\rng\varphi\subseteq X\cup\rng\varphi$,
and the inclusion map $j\colon K\into P$ is an embedding from
$\seq{K,\lambda}$ into $\seq{P,\mu}$.
By assumption on $\seq{L,\varphi}$, there exists a homomorphism
$g\colon\seq{P,\mu}\to\seq{L,\varphi}$ such that $g\circ j=f$.
Put $x=g(t)$. Then $a=b\wedge x$ and $c=b\vee x$.

The proofs of (ii)--(iv) proceed in the same way, as shown in
20.2--20.7 in \cite{Wehr}. For proving the containment
$X\cap\dnw\rng\varphi\subseteq\rng\varphi$, we need to imitate the
second part of the proof of Lemma~20.7 in \cite{Wehr}. More
specifically, let $\alpha\in X$, let $o<i$ in $L$ such that
$0<\alpha<\varphi\Theta_L(o,i)$, put $K=\set{o,i}$, let
$f\colon K\into L$ be the inclusion map, and let
$\lambda=\varphi\circ\Conc f$. Furthermore, let $P=\set{o,x,i}$ be
the three-element chain, with $o<x<i$, and let $j\colon K\into P$
be the inclusion map. Endow $P$ with the \jzh\
$\mu\colon\Conc P\to S$ defined by $\mu\Theta_P(o,x)=\alpha$ and
$\mu\Theta_P(x,i)=\varphi\Theta_L(o,i)$. Observe that the range of
$\mu$ is contained into $X\cup\rng\varphi$. By assumption on
$\seq{L,\varphi}$, there exists a homomorphism
$g\colon\seq{P,\mu}\to\seq{L,\varphi}$ such that $g\circ j=f$. Hence
the element
$\alpha=\mu\Theta_P(o,x)=(\varphi\circ\Conc g)\Theta_P(o,x)$ belongs
to the range of $\varphi$.

The proof of (iii) goes along similar lines, although the lattice $K$
and the partial lattice $P$ to be considered are much more
complicated, see  20.2--20.6 in \cite{Wehr} for details.
\end{proof}

Now we are coming to the main result (stated in the
Introduction) of Section~\ref{S:LockBr}:

\begin{proof}[Proof of Theorem~\textup{\ref{T:kappaRepr}}]
We first deal separately with the case where $\kappa=\aleph_0$,
\emph{i.e.}, $S$ is countable. Then, by Bergman's Theorem
\cite{Berg86,GoWe} and Corollary~7.5 in \cite{GoWe}, there exists a
relatively complemented \emph{modular} lattice $L$ with zero such that
$\Conc L\cong S$, moreover, if $S$ is bounded, then $L$ is bounded.

Suppose now that $\kappa>\aleph_0$. We can decompose $S$ as
$S=\bigcup_{\xi<\kappa}S_\xi$, for an increasing family
$(S_\xi)_{\xi<\kappa}$ of infinite \jz-subsemilattices of
$S$ such that $|\xi|\leq|S_\xi|<\kappa$ for all $\xi<\kappa$. Furthermore,
if $S$ is bounded, then we may assume that $1\in S_\xi$ for all
$\xi<\kappa$.

Now we construct \ML{S}s $\seq{L_\xi,\varphi_\xi}$, for $\xi<\kappa$, as
follows. For $\xi<\kappa$, suppose that $\seq{L_\eta,\varphi_\eta}$ has
been constructed for all $\eta<\xi$, such that
$\seq{L_\zeta,\varphi_\zeta}$ is an extension of
$\seq{L_\eta,\varphi_\eta}$, $|L_\eta|\leq|S_\eta|$, and
$\seq{L_\eta,\varphi_\eta}$ is lower $S_\eta$-saturated,  for
$\eta\leq\zeta<\xi$. Put $L'_\xi=\bigcup_{\eta<\xi}L_\eta$
and $\varphi'_\xi=\bigcup_{\eta<\xi}\varphi_\eta$.
Observe that $|L'_\xi|\leq|S_\xi|$. Hence, by 
Proposition~\ref{P:EmbSat} applied to $\seq{L'_\xi,\varphi'_\xi}$, there
exists a lower $S_\xi$-saturated $\seq{L_\xi,\varphi_\xi}$
with $|L_\xi|\leq|S_\xi|$ such that
$\seq{L'_\xi,\varphi'_\xi}$ admits a $0$-lattice embedding into
$\seq{L_\xi,\varphi_\xi}$ (the embedding condition is vacuously satisfied
for $\xi=0$). In particular, $L_\xi$ is a lattice with zero.
Furthermore, if $S$ is bounded, then this embedding may be taken
internal, with $L_0$ bounded and $1\in\rng\varphi_0$.

Take $L=\bigcup_{\xi<\kappa}L_\xi$, a lattice with zero. Then
$\varphi=\bigcup_{\xi<\kappa}\varphi_\xi$ is, by
Proposition~\ref{P:SummSat}, an isomorphism from $\Conc L$ onto $S$.
If $S$ is bounded, then so is~$L$.
Furthermore, by Proposition~\ref{P:SummSat}, $L$ is relatively
complemented.
\end{proof}

We observe that the lattice $L$ constructed in the proof of
Theorem~\ref{T:kappaRepr} satisfies many other properties than being
relatively complemented, such as item (iii) in the statement of
Proposition~\ref{P:SummSat}.

\section{The spaces $P^*_\kl$, $P_\kl$, $A_\kl$, $U_\kl$, $V_\kl$}
\label{S:Exmples}

For a partially ordered set $P$, we denote by $\Int P$ the Boolean
subalgebra of the powerset algebra of $P$ generated by all lower subsets of
$P$. For a limit ordinal $\lambda$, we define 
a subset $x$ of $\lambda$ to be \emph{bounded}, if $x\subseteq\alpha$ for
some $\alpha<\lambda$, and then we define a map
$\chi_\lambda\colon\Int\lambda\to\two$ by the rule
 \[
 \chi_\lambda(x)=\begin{cases}
 0,&\text{ if }x\text{ is bounded},\\
 1,&\text{ otherwise}.
 \end{cases}
 \]
We leave to the reader the easy proof of the following lemma:

\begin{lemma}
The map $\chi_\lambda$ is a $\seq{\vee,\wedge,0,1}$-homomorphism from
$\Int\lambda$ onto $\two$, for any limit ordinal $\lambda$.
\end{lemma}

For the remainder of this section, we shall fix infinite cardinals $\kappa$
and $\lambda$. Then we put
 \begin{align*}
 A_\kl&=\Int\kappa\times\Int\lambda\times\Int\lambda,\\
 U_\kl&=\setm{\seq{x_0,x_1,x_2}\in A_\kl}
 {\chi_\kappa(x_0)=\chi_\lambda(x_1)=\chi_\lambda(x_2)},\\
 V_\kl&=\setm{\seq{x_0,x_1,x_2}\in A_\kl}
 {\chi_\kappa(x_0)=\chi_\lambda(x_1)},\\
 P^*_\kl&=\setm{\seq{x_0,x_1,x_2}\in A_\kl}
 {\text{either }\chi_\kappa(x_0)=0\text{ or }x_1\cup x_2\neq\es},\\
 P_\kl&=\setm{\seq{x_0,x_1,x_2}\in A_\kl}
 {\chi_\kappa(x_0)=\chi_\lambda(x_1)\vee\chi_\lambda(x_2)}.
 \end{align*}
We endow each of the sets $A_\kl$, $U_\kl$, $V_\kl$, $P^*_\kl$, $P_\kl$
with the structure of partial lattice inherited from the (Boolean) lattice
structure of $A_\kl$, \emph{i.e.}, for a nonempty finite subset $X$ of
$P_\kl$ and $a\in P_\kl$, $a=\bigvee X$ if $a$ is the join of $X$ in
$A_\kl$, and similarly for the meet.

The following easy lemma summarizes the elementary properties of these
objects:

\begin{lemma}\label{L:AUVSP}\hfill
\begin{enumerate}
\item $A_\kl$, $U_\kl$, and $V_\kl$ are Boolean algebras such that
$U_\kl\subset V_\kl\subset A_\kl$.

\item $P^*_\kl$ is a $\seq{\vee,0,1}$-subsemilattice of
$A_\kl$ and it contains $V_\kl$.

\item For all $x$, $y\in P_\kl$, $x\sd y$ belongs to $P^*_\kl$.
\end{enumerate}
\end{lemma}

Of course, $x\sd y$ is an abbreviation for $x\wedge\neg y$.

Now let $S$ be a \jzs, let $\vec\xa=(\xa_\xi)_{\xi<\kappa}$ (resp.
$\vec\xb=(\xb_\eta)_{\eta<\lambda}$) be an increasing (resp., decreasing)
$\kappa$-sequence (resp., $\lambda$-sequence) of elements of $S$ such
that $\vec\xa\leq\vec\xb$, \emph{i.e.}, $\xa_\xi\leq\xb_\eta$ for all
$\xi<\kappa$ and all $\eta<\lambda$. We suppose, in addition, that
$\xa_0=0$.

We define a map $\sigma_\ab\colon P^*_\kl\to S$ by
the rule
 \[
 \sigma_\ab(\seq{x_0,x_1,x_2})=\begin{cases}
 \xa_{\sup x_0},&\text{ if }x_1\cup x_2=\es,\\
 \xb_{\min(x_1\cup x_2)},&\text{ otherwise}.
 \end{cases}
 \]

\begin{lemma}\label{L:sigadd}
The map $\sigma_\ab$ is a \jzh\ from $P^*_\kl$ to
$S$.
\end{lemma}

Now we define a map $\mu_\ab\colon P_\kl\times P_\kl\to S$ by the rule
 \[
 \mu_\ab(x,y)=\sigma_\ab(x\sd y),\qquad\text{for all }x,\,y\in P_\kl.
 \]
This definition is consistent, by Lemma~\ref{L:AUVSP}(iii).

\begin{lemma}\label{L:muabmeas}
The map $\mu_\ab$ is a measure (see Definition~\textup{\ref{D:meas}}) on
$P_\kl$.
\end{lemma}

\begin{proof}
This follows immediately from Lemma~\ref{L:sigadd}.
\end{proof}

By Lemma~\ref{L:meas}, there exists a unique \jzh\linebreak
$\varphi_\ab\colon\Conc P_\kl\to S$ such that
$\varphi_\ab\Theta_{P_\kl}^+(x,y)=\mu_\ab(x,y)$, for all $x$,\linebreak
$y\in P_\kl$.

Now we come to the main result of this section.

\begin{proposition}\label{P:gfablift}
Suppose that the map $\varphi_\ab\colon\Conc P_\kl\to S$ can be factored
through a lattice. Then there exists $\xc\in S$ such that
$\xa_\xi\leq\xc\leq\xb_\eta$, for all $\xi<\kappa$ and all $\eta<\lambda$.
\end{proposition}

\begin{proof}
Suppose that there are a lattice $L$, a homomorphism
$f\colon P_{\kappa,\lambda}\to L$ of
partial lattices, and a \jzh\ $\psi\colon\Conc L\to S$ such that
$\varphi_\ab=\psi\circ\Conc f$. We put
 \[
 \xc=\psi\Theta_L(f(\seq{\es,\es,\es}),
 f(\seq{\kappa,\lambda,\es})\wedge f(\seq{\kappa,\es,\lambda})).
 \]
We prove that $\xc$ satisfies the required inequalities.

Let $\xi<\kappa$. {}From the inequality
 \[
 f(\seq{\kappa,\lambda,\es})\wedge f(\seq{\kappa,\es,\lambda})
\geq f(\seq{\xi+1,\es,\es})
 \]
follows that
 \begin{multline*}
 \xc\geq\psi\Theta_L(f(\seq{\es,\es,\es}),f(\seq{\xi+1,\es,\es}))
 =\varphi_\ab\Theta_P(\seq{\es,\es,\es},\seq{\xi+1,\es,\es})\\
 =\mu_\ab(\seq{\xi+1,\es,\es},\seq{\es,\es,\es})
 =\sigma_\ab(\seq{\xi+1,\es,\es})=\xa_\xi.
 \end{multline*}
Now let $\eta<\lambda$. We first observe that
$f(\seq{\kappa,\es,\lambda})\leq f(\seq{\kappa,\lambda\sd\eta,\lambda})$ and
that $\seq{\kappa,\lambda,\es}\wedge\seq{\kappa,\lambda\sd\eta,\lambda}$ is
defined in $P_\kl$, with value $\seq{\kappa,\lambda\sd\eta,\es}$. It follows
that
 \[
 \xc\leq\psi\Theta_L(f(\seq{\es,\es,\es}),f(\seq{\kappa,\lambda\sd\eta,\es}))
=\sigma_\ab(\seq{\kappa,\lambda\sd\eta,\es})=\xb_{\min(\lambda\sd\eta)}=
\xb_\eta,
 \]
which concludes the proof.
\end{proof}

\begin{definition}\label{D:klip}
Let $P$ be a partially ordered set, let $\kappa$ and $\lambda$ be infinite
cardinals. We say that $P$ has the \emph{\klip}, if for every increasing
$\kappa$-chain $\vec\xa$ and every decreasing $\lambda$-chain $\vec\xb$ of
$P$ such that $\vec\xa\leq\vec\xb$, there exists $\xc\in P$ such that
$\vec\xa\leq\xc\leq\vec\xb$.
\end{definition}

Observe that if $\kappa=\lambda=\aleph_0$, then $P_\kl=P_{\go,\go}$ is
countable, and we obtain the following result:

\begin{proposition}\label{P:ctblecpl}
Let $S$ be a \jzs\ that does not satisfy the
$\seq{\go,\go}$-interpolation property. Then there exists a \jzh\
$\varphi\colon\Conc P_{\go,\go}\to S$ that cannot be factored through
a lattice.
\end{proposition}

\section{Necessity of the conditional completeness}\label{S:Cpl}

\begin{definition}\label{D:CondCpl}
Let $P$ be a partially ordered set. We say that $P$ is \emph{conditionally
complete}, if every nonempty majorized subset of $P$ has a least upper
bound.
\end{definition}

We recall the following elementary fact about conditional completeness:

\begin{lemma}\label{L:IPtoCC}
For any lattice $S$, if $S$ has the \klip\ for all infinite cardinals
$\kappa$ and $\lambda$, then $S$ is conditionally complete.
\end{lemma}

Then we get immediately the following result:

\begin{proposition}\label{P:Th1toCpl}
Let $S$ be a \jzs\ such that for every partial lattice $P$, every \jzh\
$\varphi\colon\Conc P\to S$ can be lifted. Then $S$ is a
conditionally complete lattice.
\end{proposition}

\begin{proof}
By \cite{TuWe2}, the condition above, even restricted to $P$ a Boolean
lattice, is sufficient to imply that $S$ is a lattice. The conclusion
follows then from Lemma~\ref{L:IPtoCC} and Proposition~\ref{P:gfablift}.
\end{proof}

In order to be able to formulate the forthcoming
Proposition~\ref{P:2amklip}, we introduce some additional notation. For
infinite cardinals $\kappa$ and $\lambda$, let
$e_\kl\colon U_\kl\into V_\kl$ be the inclusion map, let
$s_\kl\colon U_\kl\to U_\kl$, $\seq{x_0,x_1,x_2}\mapsto\seq{x_0,x_2,x_1}$ be
the natural symmetry, and put $e'_\kl=e_\kl\circ s_\kl$.

\begin{proposition}\label{P:2amklip}
Let $S$ be a \jzs, let $\kappa$ and $\lambda$ be infinite cardinal
numbers, let $\vec\xa$ (resp., $\vec\xb$) be an increasing (resp.,
decreasing) $\kappa$-sequence (resp., $\lambda$-sequence) of elements of
$S$ such that $\vec\xa\leq\vec\xb$. We denote by $\mu$ (resp., $\nu$)
the restriction of $\sigma_\ab$ to $U_\kl$ (resp., $V_\kl$).
Suppose that there are a meet-semilattice $L$, meet-homomorphisms
$f$, $f'\colon V_\kl\to L$, and an order-preserving map $\rho\colon L\to S$
such that $f\circ e_\kl=f'\circ e'_\kl$ and
$\rho\circ f=\rho\circ f'=\nu$.
Then there exists $\xc\in S$ such that $\xa_\xi\leq\xc\leq\xb_\eta$, for
all $\xi<\kappa$ and all $\eta<\lambda$.
\end{proposition}

The statement of Proposition~\ref{P:2amklip} means that if the
amalgamation problem described by the diagram below can be solved,
 \[
{
\xymatrixrowsep{2pc}\xymatrixcolsep{1pc}
\def\labelstyle{\displaystyle}
\xymatrix{
 & \seq{L,\rho} &\\
 \seq{V_\kl,\nu}\ar@{-->}[ru]^{f} & & \seq{V_\kl,\nu}\ar@{-->}[lu]_{f'}\\
 & \seq{U_\kl,\mu}\ar@{_(->}[lu]^{e_\kl}\ar@{^(->}[ru]_{e'_\kl}
}
}
 \]
for a meet-semilattice $L$, meet-homomorphisms $f$, $f'\colon V\to L$, and
an order-preserving $\rho\colon L\to S$, then there exists $\xc\in S$ such
that $\vec\xa\leq\xc\leq\vec\xb$.

\begin{proof}
Suppose that $L$, $\rho$, $f$, and $f'$
are as required. We define an element $\xc$ of $S$ by
 \[
 \xc=\rho(f(\seq{\kappa,\lambda,\es})\wedge
 f'(\seq{\kappa,\lambda,\es})).
 \]
Now we prove that $\xa_\xi\leq\xc$, for all $\xi<\kappa$. Indeed, from the
inequalities
 \begin{align*}
 f(\seq{\kappa,\lambda,\es})&\geq f(\seq{\xi+1,\es,\es})\\
 f'(\seq{\kappa,\lambda,\es})&\geq
 f'(\seq{\xi+1,\es,\es})=f(\seq{\xi+1,\es,\es})
 \end{align*}
follows that
 \[
 \xc\geq\rho(f(\seq{\xi+1,\es,\es}))=\nu(\seq{\xi+1,\es,\es})=\xa_\xi.
 \]
Next, we prove that $\xc\leq\xb_\eta$, for all $\eta<\lambda$.
Indeed, put $g=f\circ e=f'\circ e'$. We compute:
 \begin{align*}
 f'(\seq{\kappa,\lambda,\es})&\leq
 f'(\seq{\kappa,\lambda,\lambda\sd\eta})&&
 (\text{now, we observe that }
 \seq{\kappa,\lambda,\lambda\sd\eta}\in U_\kl)\\
 &=f'\circ s_\kl(\seq{\kappa,\lambda\sd\eta,\lambda})\\
 &=f(\seq{\kappa,\lambda\sd\eta,\lambda}),
 \end{align*}
thus, since $f$ is a meet-homomorphism,
 \begin{align*}
 f(\seq{\kappa,\lambda,\es})\wedge f'(\seq{\kappa,\lambda,\es})&\leq
 f(\seq{\kappa,\lambda,\es}\wedge\seq{\kappa,\lambda\sd\eta,\lambda})\\
 &=f(\seq{\kappa,\lambda\sd\eta,\es}).
 \end{align*}
Therefore,
 \[
 \xc\leq\rho(f(\seq{\kappa,\lambda\sd\eta,\es}))
 =\nu(\seq{\kappa,\lambda\sd\eta,\es})=\xb_\eta,
 \]
which completes the proof.
\end{proof}

As a corollary, even a weak version of Theorem~\ref{T:2Main} is sufficient
to require the assumption that $S$ is a conditionally complete
distributive lattice:

\begin{corollary}\label{C:CondCpl}
Let $S$ be a \jzs\ such that for every truncated square $\DD$ of
Boolean lattices and homomorphisms of Boolean lattices, every
homomorphism
$\varphi\colon\Conc\DD\to S$ can be lifted through a lattice. Then $S$
is a conditionally complete distributive lattice.
\end{corollary}

\begin{proof}
The one-dimensional version of Theorem~\ref{T:2Main} asks whether every
\jzh\ from $\Conc K$ to $S$ can be lifted, for any lattice $K$. By
\cite{TuWe2}, even the restriction of this result to the case where $K$ is
Boolean is already sufficient to imply that $S$ is a distributive
lattice. Of course, the two-dimensional amalgamation property above is
stronger (take $B_0=B_1=B_2$ and $e_1=e_2=\id_{B_0}$).

Now, if $S$ is not conditionally complete, then, by
Lemma~\ref{L:IPtoCC}, there are infinite cardinals $\kappa$ and
$\lambda$, an increasing
$\kappa$-chain $\vec\xa$ of $S$, and a decreasing $\lambda$-chain
$\vec\xb$ of $S$ such that $\vec\xa\leq\vec\xb$ but there exists no
$\xc\in S$ such that $\vec\xa\leq\xc\leq\vec\xb$. Let $\mu$ and $\nu$
be the restrictions of
$\sigma_\ab$ to $U_\kl$ and $V_\kl$, respectively, let $\varphi$ denote
the canonical isomorphism from $\Conc V_\kl$ onto $V_\kl$, and put
$\ol{\nu}=\nu\circ\varphi$. Hence $\ol{\nu}$ is a \jzh\ from $\Conc V_\kl$ to
$S$. By assumption on $S$, there are a lattice $L$, lattice
homomorphisms
$f$, $f'\colon V_\kl\to L$, and an isomorphism $\eps\colon\Conc L\to S$
such that $f\circ e_\kl=f'\circ e'_\kl$ (let us denote this map by $h$)
and $\ol{\nu}=\eps\circ\Conc f=\eps\circ\Conc f'$. Define a map
$\rho\colon L\to S$ by the rule
 \[
 \rho(x)=\eps\Theta_L(h(0_{U_\kl}),x),\qquad\text{for all }x\in L.
 \]
Then $\rho$ is order-preserving and $\rho\circ f=\rho\circ f'=\nu$, a
contradiction.
\end{proof}

\section{Lifting truncated cubes}\label{S:TruncCube}

The question whether the results of this paper can be extended from
truncated squares to truncated \emph{cubes} of lattices has a trivial,
negative answer. Indeed, let us consider the following diagram $\DD$ of
lattices and $0$-preserving lattice embeddings,
 \[
{
\def\labelstyle{\displaystyle}
\xymatrixcolsep{1pc}
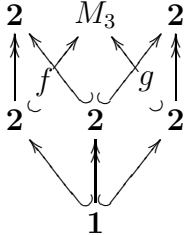
\xymatrix{
\two & M_3 & \two \\
\two\ar@{->>}[u]^{}\ar@{_(->}[ru]|-(.3){f} &
\two\ar@{_(->}[lu]^{}\ar@{^(->}[ru]^{} &
\two\ar@{->>}[u]_{}\ar@{^(->}[lu]|-(.3){g}\\
& \one\ar@{_(->}[lu]^{}\ar@{->>}[u]^{}\ar@{^(->}[ru]_{} &
}
}
 \]
where $\one=\set{0}$, $M_3=\set{0,a,b,c,1}$ is the five element modular
nondistributive lattice (with atoms $a$, $b$, $c$), the unlabeled arrows
are uniquely determined, $f(1)=a$, and $g(1)=c$. Then the image of $\DD$
under the $\Conc$ functor is obtained by truncating the top $\two$ from the
following commutative diagram of \jzs s and \jz-embeddings,
 \[
{
\def\labelstyle{\displaystyle}
\xymatrixcolsep{1pc}
\xymatrixrowsep{1pc}
\xymatrix{
& \two & \\
\two\ar[ur]_{} & \two\ar[u]_{} & \two\ar[ul]_{} \\
\two\ar[u]^{}\ar[ru]_{} &
\two\ar[lu]^{}\ar[ru]^{} &
\two\ar[u]_{}\ar[lu]_{}\\
& \one\ar[lu]^{}\ar[u]^{}\ar[ru]_{} &
}
}
 \]
that defines a homomorphism $\varphi\colon\Conc\DD\to\two$. Suppose that
$\varphi$ can be lifted to a homomorphism from $\DD$ to some
partial lattice $P$, in particular, $P$ is \emph{simple}. Let
$u\colon\two\to P$, $w\colon M_3\to P$, and $v\colon\two\to P$ be the
homomorphisms of partial lattices that correspond to the top part of such
a lifting. Travelling through the diagram~$\DD$, we obtain
 \[
 w(a)=wf(1)=u(1)=v(1)=wg(1)=w(c),
 \]
but $\Conc w$ isolates zero, \emph{i.e.}, $w$ is an embedding, hence
$a=c$, a contradiction.

Therefore, even the simplest nontrivial lattice $\two$ does not satisfy
what could be called the `three-dimensional amalgamation property'.

\section{Open problems}

The first two open problems ask whether the sufficient
conditions underlying Theorems \ref{T:1Main} and \ref{T:2Main} are also
necessary (we conjecture that yes). Possible formulations are the
following:

\begin{problem}\label{Pb:LiftPart}
Let $S$ be a distributive \jzs. If, for every partial lattice~$P$, every
\jzh\ from $\Conc P$ to $S$ can be factored through a lattice, is
$S$ \ccb?
\end{problem}

By Proposition~\ref{P:gfablift}, $S$ has the \klip, for
all infinite cardinals $\kappa$ and $\lambda$.

\begin{problem}\label{Pb:2Amalg}
Let $S$ be a distributive \jzs. If, for every truncated square $\DD$ of
lattices, every homomorphism from $\Conc\DD$ to $S$ can be factored
through a partial lattice, is $S$ \ccb?
\end{problem}

The proof of Corollary~\ref{C:CondCpl} shows that if, for every
truncated square $\DD$ of lattices, every homomorphism from $\Conc\DD$
to $S$ can be factored through a \emph{lattice}, then $S$ has the
\klip, for all infinite cardinals $\kappa$ and $\lambda$.

On the positive side, we formulate the following question, related to
Theorem~\ref{T:kappaRepr}:

\begin{problem}\label{Pb:Card}
Let $\kappa$ be an infinite cardinal, let $S$ be a \ckcb\ \jzs. Does
there exist a relatively complemented \emph{modular} lattice $L$ with
zero such that $\Conc L\cong S$?
\end{problem}

By Bergman's Theorem and the main result of \cite{Wehr00}, the answer to
Problem~\ref{Pb:Card} is known to be positive for $\kappa=\aleph_0$ and
for $\kappa=\aleph_1$.

\end{document}